\documentclass[letterpaper, 10 pt, conference]{ieeeconf}
\expandafter\let\csname ver@amsthm.sty\endcsname\relax
\let\theoremstyle\relax

\usepackage{amsthm}
\usepackage{amsmath,amsthm,epsfig,amssymb}
\usepackage{enumerate}
\usepackage{comment}
\usepackage{algorithm}
\usepackage{algorithmic}
\usepackage{amsfonts}       
\usepackage{dsfont}
\let\labelindent\relax
\usepackage{enumitem}
\usepackage{amssymb}
\usepackage[dvipsnames]{xcolor}

\IEEEoverridecommandlockouts

\DeclareMathOperator*{\argmin}{arg\,min}

\DeclareMathOperator{\dom}{dom}

\DeclareMathOperator{\intr}{int}
\DeclareMathOperator{\bdry}{bdry}

\newcommand{\bR}{\mathbb{R}}

\newcommand{\exR}{\overline{\mathbb{R}}}

\newcommand{\cC}{\mathcal{C}}

\makeatletter
\newcommand{\prox}[3][\@nil]{%
  \def\tmp{#1}%
   \ifx\tmp\@nnil
       \operatorname{prox}_{#3}^{#2}
    \else
         \operatorname{prox}_{#3}^{#1 \star #2}
    \fi}

\newcommand{\bprox}[3][\@nil]{%
  \def\tmp{#1}%
   \ifx\tmp\@nnil
       \operatorname{bprox}_{#3}^{#2}
    \else
        \operatorname{bprox}_{#3}^{#1 #2}
    \fi}
\makeatother

\usepackage{hyperref}
\hypersetup{
    colorlinks=true,
    linkcolor=NavyBlue,
    filecolor=magenta,
    citecolor =ForestGreen,  
    urlcolor=magenta,
    pdftitle={Anisotropic gradient flow}
    }
\usepackage[capitalize]{cleveref}[0.19]

\crefname{section}{section}{sections}
\crefname{subsection}{subsection}{subsections}
\Crefname{section}{Section}{Sections}
\Crefname{subsection}{Subsection}{Subsections}

\Crefname{figure}{Figure}{Figures}

\crefformat{equation}{\textup{#2(#1)#3}}
\crefrangeformat{equation}{\textup{#3(#1)#4--#5(#2)#6}}
\crefmultiformat{equation}{\textup{#2(#1)#3}}{ and \textup{#2(#1)#3}}
{, \textup{#2(#1)#3}}{, and \textup{#2(#1)#3}}
\crefrangemultiformat{equation}{\textup{#3(#1)#4--#5(#2)#6}}%
{ and \textup{#3(#1)#4--#5(#2)#6}}{, \textup{#3(#1)#4--#5(#2)#6}}{, and \textup{#3(#1)#4--#5(#2)#6}}

\Crefformat{equation}{#2Equation~\textup{(#1)}#3}
\Crefrangeformat{equation}{Equations~\textup{#3(#1)#4--#5(#2)#6}}
\Crefmultiformat{equation}{Equations~\textup{#2(#1)#3}}{ and \textup{#2(#1)#3}}
{, \textup{#2(#1)#3}}{, and \textup{#2(#1)#3}}
\Crefrangemultiformat{equation}{Equations~\textup{#3(#1)#4--#5(#2)#6}}%
{ and \textup{#3(#1)#4--#5(#2)#6}}{, \textup{#3(#1)#4--#5(#2)#6}}{, and \textup{#3(#1)#4--#5(#2)#6}}

\newtheorem{theorem}{Theorem}[section]

\newtheorem{assumption}[theorem]{Assumption}
\AddToHook{env/assumption/begin}{\crefalias{theorem}{assumption}}

\newlist{lemenum}{enumerate}{1} 
\setlist[lemenum]{label=(\roman*), ref=\theproposition(\roman*), font=\rm}
\crefalias{lemenumi}{lemma} 

\newlist{thmenum}{enumerate}{1} 
\setlist[thmenum]{label=(\roman*), ref=\theproposition(\roman*), font=\rm}
\crefalias{thmenumi}{thmenum}

\newtheorem{proposition}[theorem]{Proposition}
\newlist{propenum}{enumerate}{1} 
\setlist[propenum]{label=(\roman*), ref=\theproposition(\roman*), font=\rm}
\crefalias{propenumi}{proposition} 

\newlist{defenum}{enumerate}{1} 
\setlist[defenum]{label=(\roman*), ref=\thedefinition(\roman*), font=\rm}
\crefalias{defenumi}{definition}

\theoremstyle{remark}
\newtheorem{remark}[theorem]{Remark}

\newlist{assumenum}{enumerate}{1} 
\setlist[assumenum]{leftmargin=2.1cm,label=(A\arabic*),font=\bfseries}
\creflabelformat{assumenumi}{#2#1#3}
\crefname{assumenumi}{assumption}{assumptions}
\Crefname{assumenumi}{Assumption}{Assumptions}

\usepackage{xcolor}

\begin{document}
\title{
   Nonlinearly preconditioned gradient flows
        \thanks{
            This work was supported by the Research Foundation Flanders (FWO) PhD grant 11A8T26N and research projects G081222N, G033822N, and G0A0920N; Research Council KUL grant C14/24/103.
        }
}

\author{%
	\texorpdfstring{%
		Konstantinos Oikonomidis\thanks{
			KU Leuven, Department of Electrical Engineering ESAT-STADIUS -- %
			Kasteelpark Arenberg 10, box 2446, 3001 Leuven, Belgium
			{\sf\scriptsize
					\href{mailto:konstantinos.oikonomidis@kuleuven.be}{konstantinos.oikonomidis@kuleuven.be},%
			}%
		}
		\and
        Alexander Bodard
		\and
		Jan Quan
		\and
		Panagiotis Patrinos
	}{%
		Konstantinos Oikonomidis, Alexander Bodard, Jan Quan and Panagiotis Patrinos
	}%
}

\maketitle
\begin{abstract}
    We study a continuous-time dynamical system which arises as the limit of a broad class of nonlinearly preconditioned gradient methods. Under mild assumptions, we establish existence of global solutions and derive Lyapunov-based convergence guarantees. For convex costs, we prove a sublinear decay in a geometry induced by some reference function, and under a generalized gradient-dominance condition we obtain exponential convergence. We further uncover a duality connection with mirror descent, and use it to establish that the flow of interest solves an infinite-horizon optimal-control problem of which the value function is the Bregman divergence generated by the cost. These results clarify the structure and optimization behavior of nonlinearly preconditioned gradient flows and connect them to known continuous-time models in non-Euclidean optimization.
\end{abstract}

\section{Introduction}
Gradient flows have attracted substantial interest in both continuous-time optimization and control.
A major reason is that they offer valuable intuition and enable Lyapunov-based analyses of optimization algorithms.
A classical example is the standard gradient flow
\begin{equation} \label{eq:standard-flow}
    \dot x(t) = -\nabla f(x(t)),
\end{equation}
which can be viewed as the continuous-time limit of the gradient descent iteration
\begin{equation} \label{eq:gd}
    x^{k+1} = x^k - \gamma \nabla f(x^k),
\end{equation}
obtained by letting the stepsize \(\gamma \to 0\).

Recently, this perspective has been extended well beyond the standard gradient flow, with increasingly sophisticated dynamical systems being proposed to capture the behavior of modern optimization algorithms. 
Prominent examples include second-order dissipative systems that model methods such as Polyak’s heavy-ball algorithm and Nesterov’s accelerated gradient algorithm \cite{attouch2000heavy, su2016differential}. 
Additional developments involve \textit{mirror flows}, which arise as continuous-time analogues of mirror descent and more general Bregman-type methods \cite{krichene2015accelerated, tzen2023variational}, and \textit{proximal flows}, which correspond to the proximal point algorithm and related schemes \cite{attouch2019convergence}. 

We also draw attention to \textit{normalized flows} \cite{cortes2006finite}, which serve as continuous-time counterparts of normalized gradient methods \cite{hubler2025gradient}. 
This family of algorithms is closely connected to \textit{gradient clipping techniques} \cite{zhanggradient} and has recently attracted considerable interest due to its effectiveness in neural network training and its robustness in stochastic optimization settings with heavy-tailed noise. 

In \cite{oikonomidis2025nonlinearlypreconditionedgradientmethods}, it was shown that these methods fit within a broader family of algorithms known as \textit{nonlinearly preconditioned gradient methods} (NPGMs) \cite{oikonomidis2025nonlinearlypreconditionedgradientmethods,oikonomidis2025nonlinearly}, and also referred to as anisotropic gradient descent \cite{laude2025anisotropic} or dual space preconditioning \cite{maddison2021dual} methods.
This class is characterized by updates of the form
\begin{equation} \label{eq:npgm}
    x^{k+1} = x^k - \gamma \nabla \phi^*(\nabla f(x^k)),
\end{equation}
where $\phi: \bR^n \to \exR$ is a convex, proper and lower semicontinuous \emph{reference function}, with $\phi^*$ denoting its convex conjugate, and $f: \bR^n \to \bR$ is the continuously differentiable cost function. 
The particular choice $\phi = \tfrac{1}{2}\|\cdot\|^2$, yields \(\phi^* = \frac{1}{2} \Vert \cdot \Vert^2\) due to self-conjugacy \cite[Ex.\,11.11]{RoWe98}, such that the update rule \eqref{eq:npgm} reduces to the classical gradient descent iteration \eqref{eq:gd}.
It is straightforward to verify that gradient normalization also fits within this framework by selecting a reference function $\phi$ with bounded domain.
For instance, consider 
\(
    \phi(x) = -\varepsilon(\ln(1-\|x\|)+\|x\|),
\)
with \( \dom \phi = \{x \in \bR^n \mid \|x\| < 1\}\).
For this choice, one obtains $\nabla \phi^*(y) = \frac{y}{\|y\|+\varepsilon}$. 
Hence, substituting into \eqref{eq:npgm} yields an update that, as $\varepsilon \to 0$, approximates the classical normalized gradient method \cite[Eq.\,(NSGD)]{hubler2025gradient}. 

\begin{figure}
    \centering
    \includegraphics[width=0.8\linewidth]{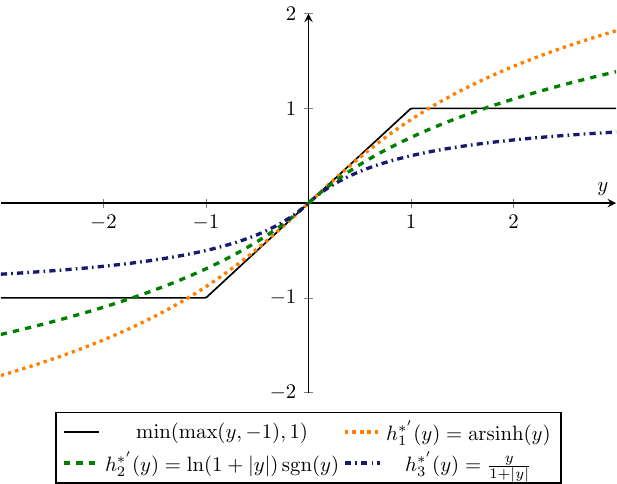}
    \caption{Preconditioners corresponding to different reference functions.}
    \label{fig:sigmoids}
\end{figure}

A variety of other reference functions \(\phi\) have been proposed.
In particular, \emph{isotropic} choices of the form \(\phi = h \circ \Vert \cdot \Vert\), for suitable scalar functions \(h\), yield algorithms that mimic gradient clipping, as illustrated in \cref{fig:sigmoids}.
A notable feature of \eqref{eq:npgm} is that its convergence is naturally analyzed under the framework of \textit{anisotropic smoothness} \cite{laude2025anisotropic}, a milder condition than standard Lipschitz smoothness that is satisfied in important applications such as matrix factorization and phase retrieval \cite{bodard2025escaping}.

In this work we address the open problem of analyzing \emph{nonlinearly preconditioned gradient flows}, defined by
\begin{equation} \label{eq:doubly_nonlinear_ode}
    \dot x(t) = -\nabla \phi^*(\nabla f(x(t))).
\end{equation}
These dynamics arise as the continuous-time analogue of the discrete-time iteration \eqref{eq:npgm} and provide a natural framework for understanding nonlinearly preconditioned gradient methods.


In fact, equations of the form \eqref{eq:doubly_nonlinear_ode} belong to a wider family of nonlinear differential equations, known as \emph{doubly nonlinear equations} \cite{stefanelli2008brezis}.
They are described by the general inclusion problem
\[
    \partial \phi(\dot u) + \partial f(u) \ni g,
\]
where in this case $u:[0,T) \to \mathcal{H}$ with $\mathcal{H}$ a real Hilbert space, $g \in L_2([0,\infty), \mathcal{H})$ and $\partial \phi$ denotes the standard subdifferential from convex analysis \cite[Def. 16.1]{bauschke2017correction}. As described in \cite{stefanelli2008brezis}, this inclusion has a straightforward mechanical interpretation: $u$ represents the displacement of a body, $g$ the corresponding energy and $\phi$ the related dissipation potential. Various works have been devoted to the study of such nonlinear dynamics, including \cite{mielke2013nonsmooth} which proposes a generalized inclusion, \cite{emmrich2011doubly} which studies a second-order system and \cite{bonaschi2016quadratic} where a stochastic model is analyzed.

\subsection{Contribution}
Our contributions can be summarized as follows:
\begin{itemize}
    \item We study the optimization properties of the nonlinearly preconditioned gradient flow \eqref{eq:doubly_nonlinear_ode} for solving smooth unconstrained minimization problems, obtaining convergence rate guarantees for suitable optimality measures. Through a novel Lyapunov-like function we retrieve standard $1/t$ rates for the suboptimality gap $f(x(t))-\inf f$ and prove the exponential convergence of the method under a generalized PL inequality.
    \item Through a duality connection with the mirror flow \cite{tzen2023variational} we obtain an optimal control viewpoint on the nonlinearly preconditioned gradient flow \eqref{eq:doubly_nonlinear_ode}. This further allows us to show that, under suitable assumptions, among all stabilizing controls, $-\nabla \phi^*(\nabla f(x))$ constitutes a control that minimizes the cost $\int_0^\infty q(x(t), u(t))$, with $q$ defined later on. The value function $V(x_0)$ that corresponds to this cost is given by the Bregman divergence generated by $f$. 
\end{itemize}

\subsection{Notation}
We denote the extended real line by $\exR = \bR \cup \{\infty\}$, the standard Euclidean inner product on $\bR^n$ by $\langle\cdot,\cdot \rangle$, and its induced norm by $\|\cdot\|$. We denote by $\mathcal{C}^k(Y)$ the class of functions that are $k$ times continuously differentiable on an open set $Y \subseteq \bR^n$. We say that a function $f:\bR^n \to \exR$ is \emph{level-bounded} if for every $\alpha \in \bR$, the set $\{x \in \bR^n \mid f(x) \leq \alpha\}$ is bounded. The \emph{effective domain} of an extended real-valued function $f:\bR^n \to \exR$ is denoted by $\dom f = \{x \in \bR^n \mid f(x)<\infty\}$. For a set $C \subseteq \bR^n$, $\intr C$ denotes its interior, the union of all open subsets contained in $C$, $\overline C$ its closure, the intersection of all closed sets on $\bR^n$ containing $C$ and $\bdry C$ its boundary, $\bdry C = \overline C \setminus \intr C$. We adopt the notions of essential smoothness, essential strict convexity and Legendre functions from \cite[Sec.\,26]{rockafellar1997convex}: we say that a proper, lsc and convex function $f:\bR^n \to \exR$ is \emph{essentially smooth} if $\intr (\dom f)\neq \emptyset$ and $f$ is differentiable on $\intr (\dom f)$ such that $\|\nabla f(x^\nu)\| \to \infty$, whenever $\intr (\dom f) \ni x^\nu \to x \in \bdry \dom f$, and \emph{essentially strictly convex}, if $f$ is strictly convex on every convex subset of $\dom \partial f$, and \emph{Legendre}, if $f$ is both essentially smooth and essentially strictly convex. In particular, a smooth convex function on $\bR^n$ is essentially smooth. For a Legendre function $f$, we denote the generated Bregman divergence $D_f(x, \bar x) = f(x)-f(\bar x) - \langle \nabla f(\bar x),x-\bar x \rangle$. We say that a function $f:\bR^n \to \exR$ is even if $f(-x)=f(x)$ for all $x \in \bR^n$.

\section{Convergence guarantees in continuous time}
In this section we study the convergence properties of the trajectory solving \eqref{eq:doubly_nonlinear_ode}, deriving sublinear convergence rates in the convex case and exponential rates under a generalized gradient dominance condition that fits the geometry of the gradient flow. We begin by formulating the optimization problem that we study
\begin{equation} \label{eq:problem}
    \min_{x \in \bR^n} f(x).
\end{equation}
Next, we present our assumptions on the problem data and our reference function.
\begin{assumption} \label{assum:main} We have the following:
    \begin{enumerate}
        \item $f \in \cC^2(\bR^n)$ is level-bounded.
        \item $\phi:\bR^n \to \exR$ is essentially smooth and strongly convex with constant $\mu_\phi > 0$. Moreover, $\phi \geq 0$ is even and $\phi(0) = 0$.
    \end{enumerate}    
\end{assumption}
We remark that due to the duality of strong convexity and Lipschitz smoothness \cite[Thm.\,18.15]{bauschke2017correction}, $\nabla \phi^*: \bR^n \to \bR^n$ is $1/\mu_\phi$-Lipschitz continuous. 
Moreover, $\nabla \phi^*$ is a bijection between $\bR^n$ and $\intr \dom \phi$ \cite[Thm.\,26.5]{rockafellar1997convex} and $\nabla \phi^*(y) = 0$ if and only if $y=0$.

Clearly, since $f \in \cC^2(\bR^n)$, $f$ is locally Lipschitz differentiable and thus the mapping $\nabla \phi^* \circ \nabla f$ is locally Lipschitz continuous. Moreover, since $f$ is level-bounded and continuous, $\argmin f \neq \emptyset$ from \cite[Thm.\,1.9]{RoWe98} and we denote $f_\star = \inf f$. 
\begin{remark}
    The assumptions on the level-boundedness of $f$ and the global Lipschitz differentiability of $\phi^*$ can be relaxed. In fact, existence of solutions can be shown under nonsmoothness of both $\phi$ and $f$ even in infinite-dimensional spaces. However, in this paper we are mostly interested in the optimization perspective, i.e., obtaining convergence rates under various conditions for $f$. For this reason, we consider meaningful practical assumptions that simplify the exposition and facilitate the analysis under standard control-theoretic arguments.
\end{remark}

We now move on to our first result that describes the global existence and uniqueness of the solution of \eqref{eq:doubly_nonlinear_ode}.
\begin{proposition}
    Let \cref{assum:main} hold true. Then, for any initial point $x(0) = x_0 \in \bR^n$, there exists a unique global solution $x:[0, \infty) \to \bR^n$ of the system \eqref{eq:doubly_nonlinear_ode}. Moreover, $\dot x\in L^2([0,\infty), \bR^n)$ and if $\phi \in \cC^2(\intr \dom \phi)$, then $\lim_{t \to \infty}\dot x(t) = 0$. 
\end{proposition}
\begin{proof}
    Since $F := \nabla \phi^* \circ \nabla f$ is locally Lipschitz, it is Lipschitz continuous on a neighborhood of $x_0$. From the Cauchy--Lipschitz (also known as Picard--Lindel\"of) theorem \cite[Thm.\,2.2]{teschl2012ordinary}, this means that there exists some $\delta > 0$ such that \eqref{eq:doubly_nonlinear_ode} admits a unique solution on $[0,\delta]$. Now note that for any $t \in [0,\delta]$ we have
    \begin{equation} \label{eq:cost_dec}
        \frac{d}{dt}(f(x(t))) = \langle \nabla f(x(t)),-\nabla \phi^*(\nabla f(x(t))) \rangle \leq 0,
    \end{equation}
    where we substituted \eqref{eq:doubly_nonlinear_ode}, and used $\nabla \phi^*(0) = 0$ and the monotonicity of $\nabla \phi^*$. Therefore, $f(x(t)) \leq f(x_0)$ for all $t\in [0,\delta]$ and since $f$ is continuous and level-bounded, $x$ is contained in a compact set. Hence, using the standard escape from compact sets argument \cite[Cor.\,2.16]{teschl2012ordinary}, we obtain the existence and uniqueness of $x$ in $[0, \infty)$.

    Now integrating \eqref{eq:cost_dec} from $0$ to $t$ we obtain
    \begin{equation*}
        \int_0^t \langle \nabla f(x(s)),\nabla \phi^*(\nabla f(x(s)) \rangle ds = f(x_0) - f(x(t)).
    \end{equation*}
    Since this equation is true for any $t \geq 0$, we have that $\int_0^\infty \langle \nabla f(x(s)),\nabla \phi^*(\nabla f(x(s)) \rangle ds \leq f(x_0) - f_\star$. Since $\phi$ is $\mu_\phi$-strongly convex, $\nabla \phi^*$ is $\mu_\phi$-cocoercive from \cite[Thm.\,18.15]{bauschke2017correction} and thus $\langle \nabla \phi^*(y),y\rangle \geq \mu_\phi \|\nabla \phi^*(y)\|^2$ for all $y \in \bR^n$. Using this inequality and $\dot x  = -\nabla \phi^*(\nabla f(x(t)))$ we then have that
    \begin{equation*}
        \int_0^\infty \|\dot x(s)\|^2 ds \leq \frac{1}{\mu_\phi}(f(x_0) - f_\star) < +\infty,
    \end{equation*}
    implying $\dot x\in L^2([0,\infty), \bR^n)$.  Assuming moreover that $\phi \in \cC^2(\intr \dom \phi)$ we have from \cite[p.\,42]{rockafellar1977higher} that $\phi^* \in \cC^2(\intr \dom \phi^*)$ and since $\dom \phi^* = \bR^n$, $\phi^* \in \cC^2(\bR^n)$. Differentiating thus \eqref{eq:doubly_nonlinear_ode}, we obtain
    \begin{align*}
        \ddot x(t) &= - \nabla^2 \phi^*(\nabla f(x(t)))\nabla^2 f(x(t))\dot x(t)
        \\
        &=  \nabla^2 \phi^*(\nabla f(x(t)))\nabla^2 f(x(t)) \nabla \phi^*(\nabla f(x(t)))
    \end{align*}
    Now note that all the functions in the right hand side are continuous and since $x$ remains in a compact set, we have that $\sup_{t \in [0,\infty)}\|\ddot x(t)\| < +\infty$, which then implies that $\dot x$ is Lipschitz continuous. From Barb{\u{a}}lat's lemma \cite[Thm.\,5]{farkas2016variations} we now obtain the claimed result.
\end{proof}
Note that in the classical setting of the gradient flow described in \eqref{eq:doubly_nonlinear_ode} where $\phi = \tfrac{1}{2}\|\cdot\|^2$, the r.h.s.\ of \eqref{eq:cost_dec} becomes $-\|\nabla f(x(t))\|^2$, thus leading to convergence guarantees for the standard stationarity measure. In our more general setting, we obtain, in light of the equality case of the Fenchel--Young inequality \cite[Prop.\,11.3]{RoWe98},
\begin{equation} \label{eq:rhs}
    \frac{d}{dt}(f(x(t))) = -[\phi^*(\nabla f(x(t))) + \phi(\nabla \phi^*(\nabla f(x(t))))],
\end{equation}
thus extending the convergence guarantees in the nonconvex setting described in \cite[Thm.\,3.2]{oikonomidis2025nonlinearlypreconditionedgradientmethods}.

We now move on to the more interesting setting where $f$ is convex. Note that obtaining convergence rates for the suboptimality gap $f(x)-f_\star$ is not straightforward and requires taking into consideration the properties of the reference function $\phi$ as well, as also noted in \cite[p.\,1006]{maddison2021dual}.
\begin{theorem} \label{thm:converge}
    Let \cref{assum:main} hold true and $x:[0,\infty)\to \bR^n$ be the unique solution of \eqref{eq:doubly_nonlinear_ode} with initial condition $x(0)= x_0$. Let, moreover, $f$ be convex. Then, we have the following:
    \begin{thmenum}
        \item $\phi^*(\nabla f(x(t)))$ is a decreasing function of $t$. \label{thm:converge:grad_dec}
        \item \label{thm:converge:gap_rate}The function $V(t) := t \phi^*(\nabla f(x(t))) + f(x(t))$ is monotonically decreasing in $[0, +\infty)$. This implies that
        \begin{equation} \label{eq:rate_conv}
            \phi^*(\nabla f(x(t))) \leq \frac{f(x_0)-f_\star}{t}, \quad \text{for} \quad t > 0.
        \end{equation}
    \end{thmenum}
    Assume moreover that $\phi = h \circ \|\cdot\|$, for some $h:\bR \to \exR$ satisfying \cref{assum:main} and let $x^\star \in \argmin f$. Then,
    \begin{thmenum}[resume]
        \item $\|\nabla f(x(t))\|$ is a decreasing function of $t$ and the same holds for $\|x(t)-x^\star\|$.\label{thm:converge:cont_fejer}
        \item If, furthermore, ${h^*}'(y)/y$ is a decreasing function on $\bR_+$, the following holds for the suboptimality gap:
         \begin{equation}
            f(x(t))-f_\star \leq \frac{\|\nabla f(x_0)\|\|x_0-x^\star\|^2}{{h^*}'(\|\nabla f(x_0)\|)t}, \text{ for } t > 0.
        \end{equation} \label{thm:converge:gap}
    \end{thmenum}
\end{theorem}
\begin{proof}
    ``\labelcref{thm:converge:grad_dec}'': We have that 
    \begin{align}
        \frac{d}{dt}(\phi^*(\nabla f(x(t)))) 
        &= \langle \nabla \phi^*(\nabla f(x(t))),\nabla^2 f(x(t))\dot x(t) \rangle \nonumber
        \\
        &\overset{\eqref{eq:doubly_nonlinear_ode}}{=} -\langle \dot x(t),\nabla^2 f(x(t)) \dot x(t) \rangle \leq 0, \label{eq:whatever-i-dont-care}
    \end{align}
    since $f$ is convex and thus $\nabla^2 f \succeq 0$.

    ``\labelcref{thm:converge:gap_rate}'': Note that 
    \begin{align*}
        \frac{d}{dt}V(t) = \phi^*(\nabla f(x(t))) + t \frac{d}{dt} &(\phi^*(\nabla f(x(t))))\\
        &+ \langle \nabla f(x(t)),\dot x(t) \rangle
    \end{align*}
    Due to \eqref{eq:whatever-i-dont-care},
    \begin{align*}
        \frac{d}{dt}V(t) &\leq
        \phi^*(\nabla f(x(t))) - \langle \nabla f(x(t)), \nabla \phi^*(\nabla f(x(t)))\rangle
        \\
        &= -\phi(\nabla \phi^*(\nabla f(x(t))))
    \end{align*}
    where we have used \eqref{eq:doubly_nonlinear_ode} and the equality case of the Fenchel--Young inequality \cite[Prop.\,11.3]{RoWe98}. Therefore, $V(t) \leq V(0)$, which means that
    \begin{equation}
        \phi^*(\nabla f(x(t))) \leq \frac{f(x_0)-f_\star}{t},
    \end{equation}
    for all $t > 0$.    

    ``\labelcref{thm:converge:cont_fejer}'': In light of \cite[Lem.\,1.3]{oikonomidis2025nonlinearlypreconditionedgradientmethods}, we have that $h^*$ is an increasing function on $\bR_+$ and $h^*(y) \geq 0$ for $y \in \bR_+$, while $\phi^* = h^* \circ \|\cdot\|$, and 
    \[
        \nabla \phi^*(y) = \frac{{h^*}'(\|y\|)}{\|y\|}y,
    \]
    for all $y \in \bR^n \setminus \{0\}$ and $0$ otherwise. Now, from \labelcref{thm:converge:grad_dec} we have that for all $t_2 \geq t_1 \geq 0$, $h^*(\|\nabla f(x(t_2))\|) \leq h^*(\|\nabla f(x(t_1))\|)$ and since $h^*$ is increasing we obtain the claimed result. For the second claim,
    \begin{align} \label{eq:dist_bound}
        \frac{d}{dt}(\tfrac{1}{2}\|x(t)&-x^\star\|^2) \nonumber
        \\
        &= \langle x(t)-x^\star, \dot x(t)\rangle \nonumber
        \\
        &= \langle x(t)-x^\star, -\nabla \phi^*(\nabla f(x(t)))\rangle \nonumber
        \\
        &= -\tfrac{{h^*}'(\|\nabla f(x(t))\|)}{\|\nabla f(x(t))\|}\langle x(t)-x^\star,\nabla f(x(t))\rangle \nonumber
        \\
        & \leq  \tfrac{{h^*}'(\|\nabla f(x(t))\|)}{\|\nabla f(x(t))\|}(f_\star - f(x(t))) \leq 0,
    \end{align}
    where the first inequality follows by the convex gradient inequality for $f$ and the last one by the fact that $f_\star \leq f(\bar x)$ for all $\bar x \in \bR^n$.

    ``\labelcref{thm:converge:gap}'': Let $W(t) = t \tfrac{{h^*}'(\|\nabla f(x_0)\|)}{\|\nabla f(x_0)\|} (f(x(t))-f_\star) + \tfrac{1}{2}\|x(t)-x^\star\|^2.$ Using \eqref{eq:dist_bound} and \eqref{eq:doubly_nonlinear_ode}, we have that
    \begin{align*}
        \frac{d}{dt}W(t) \leq -t \tfrac{{h^*}'(\|\nabla f(x_0)\|)}{\|\nabla f(x_0)\|} \langle \nabla f(x(t)),\nabla \phi^*(\nabla f(x(t))) \rangle
        \\
        - \left (\tfrac{{h^*}'(\|\nabla f(x(t))\|)}{\|\nabla f(x(t))\|} -\tfrac{{h^*}'(\|\nabla f(x_0)\|)}{\|\nabla f(x_0)\|}\right )(f(x(t)) - f_\star).
    \end{align*}
    Further, $\|\nabla f(x(t))\| \leq \|\nabla f(x_0)\|$ by \labelcref{thm:converge:cont_fejer} and since ${h^*}'(y)/y$ is decreasing, $\tfrac{{h^*}'(\|\nabla f(x(t))\|)}{\|\nabla f(x(t))\|} \geq\tfrac{{h^*}'(\|\nabla f(x_0)\|)}{\|\nabla f(x_0)\|}$. Therefore, since $f(x(t)) \geq f_\star$, $W(t) \leq W(0)$ and the claimed result follows.
\end{proof}
\begin{remark}
    The conditions imposed on $\phi$ in order to obtain convergence rates for $f(x) - f_\star$ are in fact mild: the assumption $\phi = h \circ \|\cdot \|$ covers a wide variety of methods in the related literature and is often used in order to obtain convergence guarantees. The assumption that ${h^*}'(y)/y$ is a decreasing function on $\bR_+$ actually covers a plethora of interesting reference functions as stressed in \cite{oikonomidis2025nonlinearlypreconditionedgradientmethods}. When $h \in \cC^2(\intr \dom h)$, it is equivalent to $y {h^*}''(y) \leq {h^*}'(y)$ for all $y > 0$.
\end{remark}
The standard gradient flow is known to converge exponentially in function values when $f$ satisfies a growth condition known as the PL inequality: $\|\nabla f(x)\|^2 \geq 2\mu(f(x)-f_\star)$ for all $x \in \bR^n$ and some $\mu > 0$. In our general setting, where $f$ might be nonconvex, the equivalent condition utilized to prove linear rates of convergence in the function values for the sequence generated by \eqref{eq:npgm} is known as anisotropic gradient dominance \cite[Def.\,5.6]{laude2025anisotropic}. It takes the following form: there exists some $\mu > 0$ such that for all $x \in \bR^n$,
\begin{equation} \label{eq:aniso_dom}
    \phi(\nabla \phi^*(\nabla f(x))) \geq \mu(f(x)- f_\star).
\end{equation}
Clearly, choosing $\phi = \tfrac{1}{2}\|\cdot\|^2$ we obtain the aforementioned PL inequality. Under this condition, we can show exponential convergence, as described in the following proposition.
\begin{proposition}
    Let \cref{assum:main} hold true and $x:[0,\infty)\to \bR^n$ be the unique solution of \eqref{eq:doubly_nonlinear_ode} with initial point $x_0$. Let, moreover, \eqref{eq:aniso_dom} hold. Then,
    \begin{equation*}
        f(x(t)) - f_\star \leq \exp(-\mu t)(f(x_0)-f_\star).
    \end{equation*}
\end{proposition}
\begin{proof}
    Combining \eqref{eq:rhs} and \eqref{eq:aniso_dom} and the fact that $\phi^* \geq 0$ we have that
    \begin{equation*}
        \frac{d}{dt}(f(x(t))) \leq -\mu (f(x(t))-f_\star).
    \end{equation*}
    The claimed result now follows from Gr\"onwall's inequality.
\end{proof}

\begin{remark}
    Note that the normalized gradient flow studied in \cite{cortes2006finite}, $\dot x(t)=-\tfrac{\nabla f(x(t))}{\|\nabla f(x(t))\|}$ can be brought into the generalized form of \eqref{eq:doubly_nonlinear_ode}, where equality is replaced by inclusion,
    \[
        \dot x(t) \in -\partial \phi^*(\nabla f(x(t)))
    \]
     with \(\phi^* = \|\cdot\|\) \cite[Ex.\,3.34]{beck2017first}. In this case, by \cite[Sec.\,4.4.12]{beck2017first},  $\phi$ is the indicator of the unit norm ball.
\end{remark}

\section{Mirror descent and mirror flow}
Throughout the remainder of this section we assume the following.
\begin{assumption}
    $f$ is strictly convex and supercoercive, 
    \begin{equation*}
        \lim_{\|x\| \to \infty}f(x)/\|x\| = +\infty.
    \end{equation*}
\end{assumption}
Note that since $f$ is also assumed smooth, it is Legendre, $\nabla f: \bR^n \to \intr \dom f^*$ is bijective and the set $\argmin f$ is a singleton, i.e., $\argmin f = \{x^\star\}$. Moreover, due to the duality of supercoercivity and full domain \cite[Thm.\,11.8(d)]{RoWe98}, $f^*$ has full domain.
\begin{remark}
    Strict convexity does not in general imply supercoercivity and in fact does not even imply simple coercivity as evidenced by the example $f(x)=\exp(x)$, which is a strictly convex function but $\lim_{x\to -\infty}f(x)=0$. A standard example of strictly convex functions that are also supercoercive are strongly convex ones, something which can be seen from the strong convexity inequality. On the other hand, there exist strictly convex, supercoercive functions that are not strongly convex, with a standard example being $f(x) = \|x\|^4$. 
\end{remark}

\begin{figure}
    \centering
    \vspace{0.18cm}
    \includegraphics[width=0.9\linewidth,page=1]{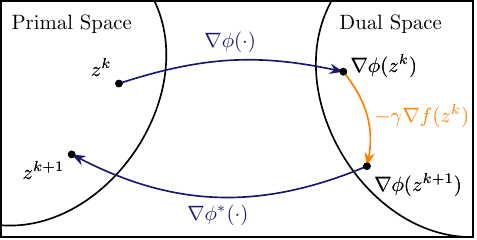}
    \caption{Visualization of a mirror descent update.}
    \label{fig:md}
\end{figure}

\begin{figure}
    \centering
    \includegraphics[width=0.9\linewidth,page=2]{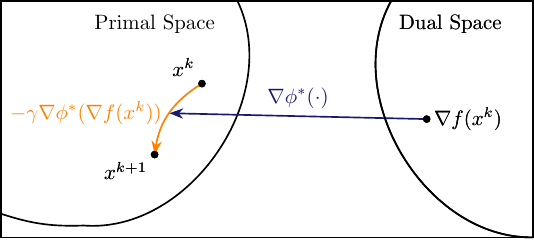}
    \caption{Visualization of a nonlinearly preconditioned gradient method update.}
    \label{fig:npgm}
\end{figure}

In discrete-time, the algorithm described in \eqref{eq:npgm} has a very interesting duality relation with the celebrated mirror descent method, which was explored in \cite{kim2023mirror}. Consider the standard mirror descent update:
\begin{equation*}
    z^{k+1} = \nabla \phi^*(\nabla \phi(z^k)-\gamma \nabla f(z^k)),
\end{equation*}
where $z^0 \in \intr \dom \phi$. Now, consider the iteration \eqref{eq:npgm} and set $z^k := \nabla f(x^k)$.
The latter is equivalent to \(x^k = \nabla f^*(z^k)\).
Substituting we obtain
\begin{equation*}
    \nabla f^*(z^{k+1}) = \nabla f^*(z^k) - \gamma \nabla \phi^*(z^k)
\end{equation*}
and thus $z^{k+1} = \nabla f(\nabla f^*(z^k) - \gamma \nabla \phi^*(z^k))$. It is thus clear that \eqref{eq:npgm} is equivalent to applying the mirror descent algorithm on $\phi^*$ with mirror potential $f$, i.e., swapping the roles of $f$ and $\phi^*$. Nevertheless, the logic behind the two methods differs. In mirror descent, the primal iterate $z^k$ is mapped to some dual iterate $\nabla \phi(z^k)$, a gradient step is performed and then it is mapped back to $z^{k+1}$ through $\nabla \phi^*$. In the dual space preconditioning framework \eqref{eq:npgm}, $\nabla f(x^k)$ is mapped to the primal space through $\nabla \phi^*$ and then a gradient step is performed directly on $x^k$ to produce $x^{k+1}$. This difference is visualized in \cref{fig:md,fig:npgm}. 

A similar relation is present also in the continuous-time analog of the two methods. To better see this, consider the gradient flow described in \cite[Eq.\,(4)]{tzen2023variational}:
\begin{equation*}
    \dot z(t) = - \nabla f(\nabla \phi^*(z(t))),
\end{equation*}
with $z(0)=z_0$. Comparing the equation above with \eqref{eq:doubly_nonlinear_ode} it is clear that one can be obtained from the other by interchanging the roles of $f$ and $\phi^*$. The described difference on the logic of the two methods is also evident in the continuous-time setup. $z(t)$ is described as the dual trajectory, while $y(t) = \nabla \phi^*(z(t))$ is the primal trajectory, the one that tends to the solutions of \eqref{eq:problem}. 

The aforementioned duality relation allows us to apply the results of \cite{tzen2023variational} directly to the setting of the current paper. To that end, let $y^\star = 0$ be the unique minimizer of $\phi^*$ (minimizer of $f$ in \cite{tzen2023variational}) and the corresponding point $x^\star = \nabla f^*(y^\star)$ ($x^\star = \nabla \phi(y^\star)$ in \cite{tzen2023variational}). In light of \cite[Thm.\,11.8]{RoWe98}, $\argmin f = \{x^\star\}$. Now consider the following control system:
\begin{equation}
    \dot x(t) = u(t),
\end{equation}
the candidate Lyapunov function $V(x) := D_f(x, x^\star)$ and $q(x,u) := \phi^*(\nabla f(x)) + \phi(-u) + \langle u,y^\star \rangle = \phi^*(\nabla f(x)) + \phi(-u)$. In the following, $q$ will play the role of the optimal control cost that we want to minimize. Now note that $V$ satisfies the properties listed in \cite[Lem.\,1]{tzen2023variational} by definition, while
\begin{align*}
    \frac{d}{dt}V(x(t)) 
    &= \langle \nabla V(x(t)),u(t) \rangle
    \\
    &= \langle \nabla f(x(t)),u(t) \rangle,
\end{align*}
and thus $\tfrac{d}{dt}V(x(t)) + q(x,u) \geq 0$ from the Fenchel--Young inequality. Therefore, denoting $J_\infty (x_0, u):=\int_0^\infty q(x(t), u(t))dt$ the cost we want to minimize over all stabilizing controls $u$ (see \cite[p.\,1543]{tzen2023variational}), i.e., functions $u:[0,\infty) \to \bR^n$ such that $\nabla f(x(t))$ converges to $0$ and thus $x(t)$ converges to $x^\star$, we retrieve the following result \cite[Thm.\,1]{tzen2023variational}:
\begin{theorem}
    For any stabilizing control $u$ and for all $t \geq 0$,
    \begin{equation*}
        \int_0^t q(x(t), u(t)) dt \geq V(x_0) - V(x(t)).
    \end{equation*}
    In particular, $J_\infty(x_0,u) \geq V(x_0)$. Moreover, for each $x_0$, the closed-loop system $\dot x(t) = -\nabla \phi^*(\nabla f(x(t)))$ gives rise to an optimal stabilizing control $u(t) = -\nabla \phi^*(\nabla f(x(t)))$, such that
    \begin{equation*}
        J_\infty(x_0, u) = V(x_0) = D_f(x_0,x^\star).
    \end{equation*}
    Furthermore, $V(x)$ is a global Lyapunov function for the closed-loop system.
\end{theorem}
Considering the simple integrator dynamics $\dot x(t) = u(t)$, the quantity $\int_0^t q(x(t), u(t)) dt \geq V(x_0) - V(x(t))$ is highly reminiscent of the energy-dissipation functional from doubly nonlinear equations, \cite[Eq.\,0.2]{pinzi2026direct}:
\begin{equation*}
    f(x(t)) + \int_0^t \left (\phi(\dot x(\tau)) + \phi^*(-\nabla f(x(\tau))) \right )d\tau - f(x_0),
\end{equation*}
highlighting once again the connection with this literature.
\begin{remark}
    The duality relation that is discussed throughout this section, allows us to obtain some convergence results directly from the analysis of the mirror descent method in~\cite{tzen2023variational}. To begin with, note that since $D_{f^*}(\nabla f(x^\star), \nabla f(x_0)) = D_f(x_0,x^\star) = f(x_0)-f_\star$, by swapping $f$ and $\phi^*$ in \cite[Thm.\,2]{tzen2023variational} (and noting that in that paper $y(t) = \nabla f(x(t))$, $y_0 = \nabla f(x_0)$), we get the result described in \labelcref{thm:converge:gap_rate}. Similarly, the second item in \cite[Thm.\,2]{tzen2023variational} describes exponential convergence under relative strong convexity of $\phi^*$ w.r.t.\ $f^*$, i.e., the following inequality holding for all $y,\bar y \in \bR^n$:
    \begin{equation*}
        D_{\phi^*}(y,\bar y) \geq \mu D_{f_\star}(y,\bar y).
    \end{equation*}
    Since this holds for all $y, \bar y\in \bR^n$, it holds also for $\bar y = \nabla f(x)$ and $y = \nabla f(x^\star) = 0$, $x \in \bR^n$. Then, $D_{f^*}(y, \bar y) = D_f(x,x^\star) = f(x) - f_\star$ and $D_{\phi^*}(y, \bar y) = \langle \nabla \phi^*(\nabla f(x)),\nabla f(x) \rangle - \phi^*(\nabla f(x)) = \phi(\nabla \phi^*(\nabla f(x)))$ and we can see that it implies \cref{eq:aniso_dom}. Nevertheless, the direct analysis of the previous section allows us to go into further detail and obtain tighter convergence guarantees under less restrictive assumptions. Note, moreover, that choosing $u =-\nabla \phi^*(\nabla f(x))$ allows us to enforce constraints on the input by choosing a suitable reference function $\phi$, i.e., one with a compatible domain. For example, $\phi(x) = \tfrac{1}{2}\|x\|^2 + \delta_{[0,1]}(\|x\|)$, leads to $u$ taking values in $B(0,1)=\{x \in \bR^n \mid \|x\| \leq 1\}$ and in this case $\nabla \phi^*$ plays the role of the projection on $B(0,1)$.
\end{remark}

\section{Conclusion}
In this paper, we analyzed nonlinearly preconditioned gradient flows, which arise as the continuous-time analogue of the the nonlinearly preconditioned gradient method. 
By means of a novel Lyapunov-like function, we established a standard sublinear convergence rate in the convex setting, and proved exponential convergence under a generalized PL inequality.
Moreover, we described a duality connection with mirror flows, which allows existing results for the mirror descent method to be directly transferred. Interesting future work involves extending our analysis to second-order equations in the spirit of \cite{attouch2000heavy} and studying the differential equation under the lens of $\Phi$-convexity similar to \cite{laude2025anisotropic}.

\scriptsize
\bibliographystyle{ieeetr}

\bibliography{references}

\end{document}